\documentclass{article}
\usepackage{amsmath,amsthm,verbatim,amssymb,amsfonts,amscd,diagrams, graphics}
\theoremstyle{plain}
\newtheorem{theorem}{Theorem}[section]

\newtheorem{lemma}[theorem]{Lemma}
\newtheorem{proposition}[theorem]{Proposition}

\theoremstyle{definition}

\theoremstyle{remark}
\newtheorem{remark}[theorem]{Remark}

\newarrow{ul}---->
\newarrow{Backwards}<----

\newcommand{\into}{\hookrightarrow}
\newcommand{\Z}{\mathbb{Z}}

\newcommand{\R}{\mathbb{R}}
\newcommand{\bd}{\partial}

\newcommand{\bp}{\boxplus}
\begin{document}

\title{All frame-spun knots are slice}
\author{Greg Friedman\\Yale University\\Dept. of Mathematics\\10 Hillhouse Ave\\PO Box 208283\\New Haven, CT 06520\\friedman@math.yale.edu\\Tel. 203-432-6473  Fax:  203-432-7316}
\date{May 13, 2003}

\maketitle

\begin{abstract}
Frame-spun knots are constructed by spinning a knot of lower dimension about a framed submanifold of $S^n$. We show that all frame-spun knots are slice (null-cobordant). 
\end{abstract}

\section{Introduction}

In \cite{L83}, Levine notes that all spun \cite{Ar} and superspun \cite{C70} knots are slice (in fact,  doubly slice) since each can be realized as a doubled disk knot. In this note, we extend the result on sliceness to the broader class of \emph{frame-spun knots} \cite{Ro89}. We use Levine's earlier criterion relating the geometric cobordism class of a higher-dimensional knot to the algebraic cobordism class of its associated Seifert matrices \cite{L69}.

By a knot, $K$, we mean a smooth (or PL-locally flat) codimension-two embedding $S^{k-2}\subset S^k$, $k\geq 3$, where $S^k$ is given the standard differential (PL) structure, but $S^{k-2}$ is only required to be \emph{homeomorphic} to the standard $(k-2)$ sphere. We will follow the usual abuse of notation and use $K$ to stand for either the embedding, its image, or their equivalence class under orientation preserving diffeomorphism of $S^k$. Two knots $K_0$, $K_1$ are cobordant if there exists a smooth (or PL) oriented manifold $V$, homeomorphic to $S^{k-2}\times I$, embedded properly and smoothly (or PL locally-flatly) in $S^{k}\times I$ such that $\bd V=(1\times K_1)\cup(0\times -K_0)$. A knot is null-cobordant, or \emph{slice}, if it is cobordant to the standard unknotted embedding of $S^{k-2}\subset S^k$; this is equivalent to the condition that $K$ is  the boundary sphere knot of a proper disk knot $D^{k-1}\subset D^{k+1}$. The set of cobordism classes of knots  $S^{k-2}\subset S^k$ forms an abelian group $C_{k-2}$ under the operation of knot connected sum.  A knot $S^{k-2}\subset S^{k}$ is doubly sliced if it is the intersection of a standard unknotted  $S^k\subset S^{k+1}$ with a \emph{trivially} knotted $S^{k-1}\subset S^{k+1}$ (in particular each resulting half of $S^{k+1}$ cut along the standard $S^{k}$ will contain a slicing disk of the knot given by the piece of $S^{k-1}$ it contains).

In \cite{Ar}, Artin introduced the technique of spinning a knot to create a knot of higher dimension. Generalizations followed, including the superspinning construction of Cappell \cite{C70} for spinning a knot about a $p$-sphere, $p\geq 1$. It was shown by Levine \cite{L83} that all spun and superspun knots are doubly slice, hence slice. In \cite{Ro89}, Roseman further generalized these spinning constructions to frame-spinning, by which knots can be spun about any closed manifold $M$ embeddable with framing inside a sphere. We show here that all frame-spun knots are also slice and provide one sufficient condition on the embedded manifold for  a knot spun about it to be doubly slice. The question of whether all frame-spun knots are doubly slice remains open; it is not hard to imagine that one could construct knots that do not satisfy Levine's necessary algebraic conditions for double sliceness \cite{L83} since the manifold $M$ no longer needs have trivial cup products, but such a construction would require being able to find a manifold/knot pair satisfying a suitable combination of algebraic properties. 

\section{Frame spinning}

We begin by describing the frame-spinning construction introduced by Roseman \cite{Ro89}; the term ``frame spinning'' was introduced by Suciu \cite{Su92}.

Suppose that $S^{m+k-2}$ is embedded in
$S^{m+k}$ by the standard (unknotted) embedding.  Let $K$ be a knot $S^{k-2}\subset S^k$, and let $M^m$ be an $m$-dimensional closed  framed submanifold of $S^{m+k-2}$ with framing $\phi$.  Note that these conditions ensure the orientability of $M$ and that the signature of $M$ will be $0$\footnote{If $M$ is smooth, this follows from the fact that $M$ must be stably parallelizable: we have $\epsilon^{m+k-1}\cong T(S^{m+k-2})\oplus \epsilon \cong T(M)\oplus \epsilon^{k-2}\oplus \epsilon$. Thus  all non-trivial characteristic classes will vanish, which implies orientability by the vanishing of the first Stiefel-Whitney class and vanishing of the signature by the Hirzebruch signature theorem. In fact, given our assumptions, the orientability and vanishing of the signature will hold even without the need to assume $M$ smooth: In the language of Fadell \cite{Fa}, our embedded tubular neighborhood pair $(\phi(M\times D^{k-2}), \phi(M\times D^{k-2})-\phi(M))$ is fiber homotopy equivalent to the \emph{associated normal fiber space} $(N, N_0)$ of the embedding. The triviality of the neighborhood  pair clearly implies its $(\Z-)$ orientability. Then, since our embedding must be locally-flat, by \cite[Thm. 4.11]{Fa}, $(T,T_0)\oplus (N, N_0)$ is fiber homotopy equivalent to $(T^*, T^*_0)$, where $T$ and $T^*$ are the respective \emph{tangent fiber spaces} of $M$ and $S^{m+k-2}$.  But by Fadell's Prop. 3.17, $(T^*, T^*_0)$ is fiber homotopy equivalent to the ordinary tangent bundle of $S^n$, which is also certainly orientable. Then by Fadell's Prop. 7.7, $(T, T_0)$ is orientable, and so by his Prop. 3.16, $M$ is orientable. The claim about the signature follows more clearly from \cite{Sc} via the theory of Pontrjagin classes for topological manifolds and the accompanying version of the signature theorem.}. 

Next, let $(D_{-}^{k}, D_{-}^{k-2})$ be an unknotted open disk pair which is
the open regular neighborhood pair of a point  of 
the knot $K\subset S^{k}$. Let $(D_{+}^{k}, D_{+}^{k-2})=(S^{k}, K)-
(D_{-}^{k},
D_{-}^{k-2})$. This is a disk knot with the unknot as its boundary. Let $M^m\times (D^k,D^{k-2})$ be the normal bundle pair of $M^m\subset S^{m+k-2}$
determined by the framing $\phi$ and the standard framing of $S^{m+k-2}$ in $S^{m+k}$. 

Now define $\sigma_M^{\phi}(K)$  to be the $(m+k-2)$-sphere 
\begin{equation*}
(S^{m+k-2}-M^m\times \text{ int } D^{k-2}) \cup_{M^m\times S^{k-3}} M^m \times D_{+}^{k-2}
\end{equation*}
embedded in the $(m+k)$-sphere
\begin{equation*}
(S^{m+k}-M^m\times \text{ int } (D^{k})) \cup_{M^m\times S^{k-1}} M^m \times D_{+}^{k}.
\end{equation*}
This construction corresponds to removing, for each point of $M$, the trivial disk pair 
$(D^{k},
D^{k-2})$, which is the fiber of the normal bundle pair of $M$, and replacing it
with the knotted disk pair 
$(D_{+}^{k}, D_{+}^{k-2})$ determined by $K$   

In the case where $M^m$ is the circle $S^1 $ with
the standard unknotted embedding and bundle framing, $\sigma^{\phi}_M(K)$ is  the Artin spin of $K$. Similarly, if $M=S^m$ with
the standard unknotted embedding and bundle framing, we obtain the
$m$-superspin of $K$  \cite{C70}. Note also that frame spinning can be further generalized to include twisting and other other deformations; see \cite{GBF}, \cite{GBF1}, or \cite{GBF7}

\section{Cobordism}

Let $K$ be a knot $S^{2n-1}\subset S^{2n+1}$. Then $K$ bounds an embedded oriented bicollared connected $2n$-manifold $V$ in $S^{2n+1}$, called a Seifert surface for $K$, and there is defined a Seifert linking pairing $F_n(V)\otimes F_n(V)\to \Z$, where $F_n(V)$ represents the homology group $H_n(V;\Z)$ modulo torsion. For chains $\alpha$ and $\beta$ in $C_n(V)$ representing such homology classes, the pairing is given by $\alpha\otimes \beta \to L(\alpha, i_*(\beta))$, where $i$ is the translation off $V$ in the positive normal direction of the bicollar and $L: F_n(V)\otimes F_n(S^n-V)\to \Z$ is the linking pairing in $S^n$. The pairing $L$ is equivalent to the Poincare-Lefschetz intersection pairing $ F_{n+1}(S^n, V)\otimes F_n(S^n-V)\to \Z$, and so at the chain level, $L(\alpha, i_*(\beta))$ can be described by the algebraic intersection number $z\cdot i(\beta)$ of representative chains in dual cell divisions such that $z\in C_{n+1}(S^{2n+1})$ and $\bd z=\alpha$. The matrix $A$ of the Seifert linking pairing is called the Seifert matrix of the knot, and it is well-defined by the knot up to $S$-equivalence (see \cite{L70}).  The matrix $A$ also enjoys the property that $A+(-1)^nA'$ is integrally unimodular, where $A'$ represents the transpose of $A$. In fact, this last matrix represents the intersection form $F_n(V)\otimes F_n(V)\to \Z$, which is non-singular since $\bd V\cong S^{2n-1}$. 

In \cite{L69}, Levine defines a square integral matrix $N$ to be null-cobordant if it is integrally congruent to a matrix of the form $
\begin{pmatrix}0 & N_1\\N_2 & N_3\end{pmatrix}$, where the $N_i$ are square matrices of the same dimensions. Two matrices $A$ and $B$ are cobordant if the block sum $A\boxplus -B$ is null-cobordant.  For fixed $\epsilon$, $\epsilon=\pm 1$, it is shown that in the set of matrices for which $N+\epsilon N'$ is integrally unimodular, cobordism is an equivalence relation, and the set of cobordism classes of matrices in this set forms an abelian group under block sum. This group is denoted $G_{\epsilon}$. Furthermore, it is shown that there is a well-defined homomorphism $\phi_n$ from $C_{2n+1}$, the group of cobordism classes of knots  $S^{2n-1}\subset S^{2n+1}$ under knot sum, to $G_{(-1)^n}$, given by assigning to the cobordism class of a knot any Seifert matrix for that knot. The main theorem of \cite{L69} is that the homomorphism $\phi_n$ is an isomorphism for $n\geq 3$, an epimorphism for $n=1$, and, for $n=2$ an isomorphism onto the subgroup of index two of $G_{1}$ consisting of  matrices such that $A+A'$ has signature a multiple of $16$.

\begin{remark}
It was shown earlier by Kervaire \cite{Ke} that $C_{2n}=0$, i.e. all even dimensional knots are cobordant to each other and hence, in particular, are cobordant to the trivial knot. In other words, all even dimensional knots are slice.
\end{remark}

\section{All frame spun knots are slice}

Suppose that $K$ is a knot $S^{k-2}\subset S^k$ and that $M^m$, is a closed $m$-manifold, $m>0$, that can be embedded with framing in $S^{2n-1}=S^{m+k-2}$. Then we can form the frame spun knot $\sigma(K)$, $S^{2n-1}\subset S^{2n+1}$, $m+k=2n+1$. Since $k\geq 3$ and $m\geq 0$, we must have $n>1$, so Levine's machinery will apply. It is shown by Klein and Suciu in \cite{ KS91} that if $V$ is a Seifert surface for $K$, then we can construct a Seifert surface $V^{\sigma}$ for $\sigma(K)$ of the form $$(D^{2n}-M\times \text{int}(D^{k-1}))\cup_{M\times D^{k-1}} M\times V.$$ The first term represents the standard ball $D^{2n}\subset S^{2n+1}$ minus the neighborhood of $M$ in its boundary, while the second term is given by spinning $V$ along with $K$. Applying the K\"unneth theorem and some obvious homotopy equivalences, the Mayer-Vietoris sequence for this union will have the form
\begin{equation*}
\to H_{n}(M)\to \underset{a+b=n}{\bigoplus} H_a(V) \otimes H_b(M) \oplus \underset{a+b=n-1}{\bigoplus} H_a(V)*H_b(M) \to  H_n( V^{\sigma}) \to .
\end{equation*}
But in each dimension, the restricted homomorphism from $H_i(M)\cong H_0(*)\otimes H_i(M)$ to the summand $H_0(V)\otimes H_i(M)$ is  an isomorphism, so quotienting out torsion, $F_n(V^{\sigma})\cong \oplus _{a+b=n, a>0} F_a(V) \otimes F_b(M)$. We will let $F_{a,b}=F_a(V) \otimes F_b(M)$. If we choose bases for each $F_a(V)$ and $F_b(M)$, we can form a basis of  $F_{a,b}$ by the elements $x \otimes \xi$, where $x\in F_a(V)$, $\xi\in F_b(M)$ are basis elements. Furthermore, we can assume that each $x\otimes \xi$ is represented by a chain of the form $x\times \xi$ (we abuse notation slightly by letting letters stand for both chains and homology classes).  If we order the bases of $F_a(V)$ and $F_b(M)$, we can order the basis for $F_{a,b}$ lexicographically.

Our goal now is to examine the Seifert matrix of $\sigma(K)$ with respect to this basis. For basis elements of the form $x\otimes \xi\in F_{a,b}$ and $y\otimes \eta\in F_{c,d}$, the corresponding matrix entry will be the linking number in $S^{2n+1}$ given by $L(x\otimes \xi, i_*(y\otimes \eta))$, which we have noted will be the intersection number of a chain bounded by $x\times \xi$ with a translate of  $y\times \eta$. If $z\in C_{a+1}(S^{k})$ is a chain with $\bd z =x$, then $z\times \xi$ will be a chain in $S^k\times M\subset S^{2n+1}$ with $\bd (z\times \xi)=x\times \xi$. Also note that $i(y\times \eta)=i(y)\times \eta$. Thus,
\begin{equation}\label{E: pair}
L(x\otimes \xi, i_*(y\otimes \eta))=(z\times \xi)\cdot (i(y)\times \eta )=(-1)^{(m-b)(k-c)}(z\cdot i(y))(\xi\cdot \eta),
\end{equation}
 where the first $\cdot $ is the intersection in $S^{k}$ and the second is the intersection number in $M$ (see \cite[VIII.13.13]{Dold}). In order for this number to be non-zero, it is necessary that $k\leq a+1 +c$ and $m\leq b+d$. Substituting $a+b=c+d=n$ and $m+k=2n+1$, we see that this intersection number must be $0$ unless $a+c=k-1$. Thus we obtain a Seifert matrix of the form
\begin{equation*}
A^{\sigma}=\begin{pmatrix}
0 &&&& L_1\\
 &&&\cdot&\\
&&\cdot\\
& \cdot\\
L_{k-2}&&&&0
\end{pmatrix},
\end{equation*}
where the $L_i$ are matrix blocks on the off-diagonal. The matrix $L_i$ is the matrix of the Seifert pairing restricted to $F_{i,n-i}\otimes  F_{k-1-i, n-(k-1-i)}$. Since $(n-i)+(n-(k-1-i))=m$, the ranks of $F_{i,n-i}$ and $F_{k-1-i, n-(k-1-i)}$ are equal, employing the Poincare dualities of $M$ and $V$ (this is acceptable for $V$ since $\bd V\cong S^{k-2}$ and the dimension ranges for $V$ are $0<i<k-1$).  Hence each matrix $L_i$ is square, and furthermore it has the same dimensions as $L_{k-1-i}$. (Note, these matrices are not necessarily transposes, although they do share a relation that becomes clear in the study of Alexander modules of knots (see \cite{L66}). However, we shall not need this here.) It is also possible that some of the blocks $L_i$ will be empty, in particular if  we  exceed the allowable
 dimension range for $M$. However, the duality relationships still hold for these empty blocks, so they occur in pairs, or possibly in the middle if $k-2$ is odd. 

The following is now immediate:

\begin{proposition}If $k$ is even (hence $m$ odd) or if $k$ is odd and $F_{\frac{k-1}{2}, \frac{k-1}{2}}$ is trivial, then $A^{\sigma}$ is null-cobordant and hence $\sigma(K)$ is slice. 
\end{proposition}

We let $[\Xi]$ stand for the cobordism class of the matrix $\Xi$.

\begin{lemma}
If $k$ is odd, $[A^{\sigma}]=[L_{\frac{k-1}{2}}]$.
\end{lemma}
\begin{proof}
Let $\mu=\sum_{i=1}^{\frac{k-3}{2}}\text{rank}(F_{i,n-i})=\sum_{i=\frac{k+1}{2}}^{k-2}\text{rank}(F_{i,n-i})$ and $\nu = \text{rank}(F_{\frac{k-1}{2}, \frac{k-1}{2}})$. Let $I_{r}$ denote the $r\times r$ identity matrix. Let $J$ be the matrix 
\begin{equation*}
\begin{pmatrix}
0 & I_{\nu} &0\\
I_{\mu} &0 &0\\
0 &0& I_{\mu} 
\end{pmatrix}.
\end{equation*}
Then
\begin{equation*}
JA^{\sigma}J'=
\begin{pmatrix}
L_{\frac{k-1}{2}} &0 &0\\
0 &0 &C\\
0 & B & 0
\end{pmatrix},
\end{equation*}
where
\begin{equation*}
B=
\begin{pmatrix}
0
 &&&& L_{\frac{k+1}{2}}\\
 &&&\cdot&\\
&&\cdot\\
& \cdot\\
L_{k-2}&&&&0
\end{pmatrix}
\qquad
C=
\begin{pmatrix}
0 &&&& L_{1}\\
 &&&\cdot&\\
&&\cdot\\
& \cdot\\
L_{\frac{k-3}{2}}&&&&0
\end{pmatrix}.
\end{equation*}
Thus $JA^{\sigma}J'$ is the block sum of $L_{\frac{k-1}{2}}$ with a null-cobordant matrix and so represents $[L_{\frac{k-1}{2}}]$ in the group.
\end{proof}

So, we turn to examine $L_{\frac{k-1}{2}}$. From equation \eqref{E: pair}, it follows that $L_{\frac{k-1}{2}}=(-1)^{\frac{mk+m}{4}} A\otimes \tau$, where $A$ is the Seifert matrix for $K$ and $\tau$ is the intersection pairing matrix on $F_{\frac{m}{2}}(M)$. 

\begin{proposition}
If $k$ is odd (hence $m$ even), then $[A^{\sigma}]=0$ in $G_{(-1)^n}$.
\end{proposition}
\begin{proof}
We already know $[A^{\sigma}]=[L_{\frac{k-1}{2}}]$, which by the above observation is $[(-1)^{\frac{mk+m}{4}} A\otimes \tau]$. If $\dim (M)\equiv 2 \mod 4$, its intersection pairing is  non-singular skew symmetric, so we are free to choose a symplectic basis for $F_{m/2}(M)$ with respect to which $\tau=\begin{pmatrix}0 & I\\ -I & 0\end{pmatrix}$. Thus, if $\text{rank}(F_{m/2}(M))=2s$ and $\bp^s A$ denotes the block sum of $s$ copies of $A$, we have 
\begin{equation*}
[(-1)^{\frac{mk+m}{4}} A\otimes \tau]=
(-1)^{\frac{mk+m}{4}}\left[\left(
\begin{matrix}0 & \bp^s A\\
- \bp^s A & 0
\end{matrix}\right)\right]=0.
\end{equation*}

Now suppose $\dim (M)\equiv 0 \mod 4$. Note that in this case the sign $(-1)^\frac{mk+m}{4}$ must be $+1$, so it can be ignored.

For any matrix $B$, the block sum $B\bp -B$ is nullcobordant; just conjugate by $J=\begin{pmatrix}I& I\\ I &  0\end{pmatrix}$.
Thus in $G_{(-1)^n}$, we have 
\begin{equation*}
[A\otimes \tau]=
[A\otimes \tau]+\left[\left(\begin{matrix}A & 0\\ 0 & -A \end{matrix}\right)\right]
=[A\otimes \tau]+ [A]+ [-A]= [A\otimes (\tau \bp 1 \bp -1)].
\end{equation*}
 But $(\tau \bp 1 \bp -1)$ is conjugate by a change of basis to a block sum $(\bp^u 1) \bp( \bp^v -1)$ (see \cite[II.4]{MH}). Since $M$ is stably parallelizable, its signature must be $0$ and so $u=v$. Therefore, $[A\otimes \tau]= [A\otimes \left((\bp^u 1) \bp (\bp^u -1)\right)]=[\bp^u A]+ [\bp^u-A]=[\bp^u A]-[\bp^u a]=0$. 
\end{proof}

Putting together the previous propositions, we obtain our main theorem:

\begin{theorem}
All frame-spun knots  are slice.
\end{theorem}

We conclude with a sufficient condition to ensure that a frame-spun knot will be doubly slice.

\begin{proposition}
Let $R:S^{n-2}\to S^{n-2}$ be the standard reflection through the  hyper-plane $x_1=0$ in $\R^{n-1}$. Suppose that $M^m$ admits a homeomorphism $r: M\to M$ and an embedding with framing $\phi: M\times D^{k-2}\into S^{n-2}$ such that $R\phi(x,y)=\phi(r(x),y)$. Then for any knot $K: S^{k-2}\subset S^k$,  $\sigma_M^{\phi}(K)$ is doubly slice. 
\end{proposition}
\begin{proof}
If the hypotheses of the theorem are met, then $\sigma_M^{\phi}(K)$ can be written as a doubled disk knot, i.e. the union along their common boundary of a disk knot $D^{n-2}\subset D^n$ and its reflection across the sphere. All doubled disk knots are doubly slice \cite{L83}.
\end{proof}

\bibliographystyle{amsplain}
\bibliography{bib}

\end{document}